\documentclass[11pt,reqno,a4paper]{amsart}

\usepackage{amsfonts, amsmath, amsthm, amssymb}
\usepackage[english]{babel}
\usepackage[mathcal]{eucal}
\usepackage{geometry}
\usepackage{relsize}
\usepackage{xcolor}
\usepackage{longtable}
\usepackage{caption}
\usepackage{setspace}
\usepackage{hyperref}

\title[]{Conditional bias reduction can be dangerous:\\
 a key example from sequential analysis}
\author{Ben Berckmoes, Anna Ivanova, Geert Molenberghs}
\keywords{conditional MLE, marginal MLE, group sequential trial, mean absolute error}
\thanks{Ben Berckmoes is post doctoral fellow at the Fund for Scientific Research of Flanders (FWO)}
\thanks{Financial support from the IAP research network \#P7/06 of the Belgian Government (Belgian Science Policy) is gratefully acknowledged.}

\date{}

\begin{document}

\maketitle

\newtheorem{pro}{Proposition}
\newtheorem{lem}[pro]{Lemma}
\newtheorem{thm}[pro]{Theorem}
\newtheorem{de}[pro]{Definition}
\newtheorem{co}[pro]{Comment}
\newtheorem{no}[pro]{Notation}
\newtheorem{vb}[pro]{Example}
\newtheorem{vbn}[pro]{Examples}
\newtheorem{gev}[pro]{Corollary}
\newtheorem{vrg}[pro]{Question}
\newtheorem{rem}[pro]{Remark}
\newtheorem{lemA}{Lemma}

\begin{abstract}
We present a key example from sequential analysis, which illustrates that conditional bias reduction can cause infinite mean absolute error.
\end{abstract}

\section{Introduction}

The following group sequential paradigm has been studied extensively in the literature, see e.g. \cite{BIM18, C89, EF90, FDL00, HP88, LH99, MKA14, W92}.

Let $X_1, X_2, \ldots$ be independent and identically distributed observations with normal law $N(\mu,\sigma^2)$, and, for each $n \in \mathbb{N}_0$, $N_n$ an $\{n,2n\}$-valued random sample size that solely depends on $X_1, \ldots, X_n$
through the stopping rule 
\begin{equation}
\mathbb{P}[N_n = n \mid X_1, \ldots, X_n] = \psi\left(K_n/n^\gamma\right),\label{eq:StoppingRule}
\end{equation}
where $\psi$ is a Borel measurable map of $\mathbb{R}$ into $[0,1]$, $\gamma \in \mathbb{R}^+_0$ a shape parameter, and, for $m \in \mathbb{N}_0$, $K_m = \sum_{k = 1}^m X_k$. The choice $\gamma = 1/2$ leads to Pocock boundaries (\cite{P77}) and the choice $\gamma = 0$ to O'Brien-Fleming boundaries (\cite{OF79}).
 
The above setting models the idea that, after having collected the data $X_1, \ldots, X_n$, it is decided, based on the stopping rule (\ref{eq:StoppingRule}), if the trial is stopped (that is, the final sample size is $N_n = n$), or continued (that is, the additional data $X_{n + 1}, \ldots, X_{2n}$ are collected and the final sample size is $N_n = 2n$). 

Assuming $\sigma$ known, the following estimators for the location parameter $\mu$ are often discussed in the literature (\cite{FDL00},\cite{MKA14}):

(a) the {\em marginal MLE}, defined by the parameter value that maximizes the marginal likelihood
\begin{equation}
\mathcal{L}(\theta;X_1,\ldots,X_{N_n}) = \frac{1}{\sigma^{N_n}}\prod_{k=1}^{N_n} \phi\left(\frac{x_k - \theta}{\sigma}\right),\label{eq:SimpleLikelihood}
\end{equation}
where $\phi$ is the standard normal density. Of course, the marginal MLE is the ordinary sample mean 
\begin{equation}
\widehat{\mu}_{N_n} = \frac{1}{N_n} \sum_{k = 1}^{N_n} X_k.\label{eq:SolutionSimpleLikelihood}
\end{equation}
This approach is simple, because it is based on the likelihood of the collected data only, without taking the stopping mechanism into account. The marginal MLE has been criticized in the literature, because it has potentially large bias (\cite{EF90}). However, it was shown in \cite{BIM18} that in many cases the bias vanishes quickly if $n$ grows.

(b) the {\em conditional MLE} $\widehat{\mu}_{c,N_n}$, defined by the parameter value that, for $N_n = m$, maximizes the conditional likelihood  
\begin{equation}
\mathcal{L}(\theta; X_1,\ldots, X_{m} \mid N_n = m) = \frac{1}{\sigma^{m}}\prod_{k=1}^{m} \phi\left(\frac{X_k - \theta}{\sigma}\right) \frac{\mathbb{P}_\theta[N = m \mid X_1,\ldots,X_m]}{\mathbb{P}_\theta[N_n = m]}.\label{eq:ComplexLikelihood}
\end{equation}
This approach is complex, because contrary to the marginal MLE, it also models the stopping mechanism. An explicit value for the conditional MLE cannot be obtained, and one has to rely on numerical methods to calculate it. However, the conditional MLE, also known as the {\em conditional bias reduction estimate} (\cite{FDL00}), is favored by the literature because it is claimed to reduce bias by taking all information into account.

In this paper, we will show that if we take $\mu = 0$, $\sigma =1$, $\psi = 1_{\left[0,\infty\right[}$, and $\gamma$ arbitrary, then
$$ \lim_{n \rightarrow \infty} \mathbb{E}[\left|\widehat{\mu}_{N_n}\right|] = 0 \text{ and } \forall n \in \mathbb{N}_0 : \mathbb{E}[\left|\widehat{\mu}_{c,N_n}\right|] = \infty.$$
That is, conditional bias reduction can cause infinite mean absolute error.

\section{Mean absolute error}

We keep the setting of the previous section, and we take $\mu = 0$, $\sigma = 1$, $\psi = 1_{\left[0,\infty\right[}$, and $\gamma$ arbitrary. So the stopping rule (\ref{eq:StoppingRule}) is now turned into 

\begin{displaymath}
\mathbb{P}[N_n = n \mid X_1, \ldots, X_n] = \left\{\begin{array}{clrr}      
1 &\textrm{ if }& K_n  \geq 0\\       
0 &\textrm{ if }& K_n < 0
\end{array}\right..
\end{displaymath}
That is, after having collected the $N(0,1)$-data $X_1,\ldots,X_n$, the trial is stopped if $K_n \geq 0$ and continued otherwise.

We first focus on the marginal MLE $\widehat{\mu}_{N_n} = K_{N_n}/N_n$. Let $\phi$ be the standard normal density and $\Phi$ the standard normal cumulative distribution function. Following \cite{BIM18}, we see that the joint density of $N_n$ and $K_{N_n}$ is given by
\begin{equation}
f_{N_n,K_{N_n}}(n,k) = \frac{1}{\sqrt{n}} \phi\left(\frac{k}{\sqrt{n}}\right) 1_{[0,\infty[}(k)\label{eq:JointDensity1}
\end{equation}
and 
\begin{eqnarray}
\lefteqn{f_{N_n,K_{N_n}}(2n,k)} \label{eq:JointDensity2}\\
&=& \frac{1}{\sqrt{2n}} \phi\left(\frac{k}{\sqrt{2n}}\right) - \int_{0}^\infty \frac{1}{n} \phi\left(\frac{u}{\sqrt{n}}\right) \phi\left(\frac{k - u}{\sqrt{n}}\right)  du\nonumber\\
&=&  \frac{1}{\sqrt{2n}} \phi\left(\frac{k}{\sqrt{2n}}\right) \left[1 - \Phi\left(\frac{k}{\sqrt{2n}}\right) \right].\nonumber
\end{eqnarray}
We learn from (\ref{eq:JointDensity1}) and (\ref{eq:JointDensity2}) that 
\begin{eqnarray}
\lefteqn{\mathbb{E}[\left|\widehat{\mu}_{N_n}\right|]}\label{eq:MAEExplicit}\\
&=& \frac{1}{n} \mathbb{E}\left[\left|K_n\right|1_{\left\{N_n = n\right\}}\right] + \frac{1}{2n} \mathbb{E}\left[\left|K_{2n}\right|1_{\left\{N_n = 2n\right\}}\right]\nonumber\\
&=& \frac{1}{n} \int_0^\infty \left|k\right| \frac{1}{\sqrt{n}} \phi\left(\frac{k}{\sqrt{n}}\right) dk + \frac{1}{2n} \int_{-\infty}^\infty \left|k\right| \frac{1}{\sqrt{2n}} \phi\left(\frac{k}{\sqrt{2n}}\right) \left[1 - \Phi\left(\frac{k}{\sqrt{2n}}\right)\right] dk\nonumber\\
&=& \frac{1}{\sqrt{n}}\left(\mathbb{E}\left[\xi1_{[0,\infty[} (\xi)\right] + \frac{1}{\sqrt{2}} \mathbb{E}\left[\left|\xi\right|\left[1 - \Phi(\xi)\right]\right]\right),\nonumber
\end{eqnarray}
with $\xi$ a standard normally distributed random variable. It clearly follows from (\ref{eq:MAEExplicit}) that 
\begin{equation}
\lim_{n \rightarrow \infty} \mathbb{E}[\left|\widehat{\mu}_{N_n}\right|] = 0.\label{eq:asymMAEmarginal}
\end{equation}
That is, the mean absolute error of $\widehat{\mu}_{N_n}$ with respect to the true parameter $0$ vanishes if $n \to \infty$.  

We now turn to the conditional MLE $\widehat{\mu}_{c,N_n}$, which maximizes the conditional likelihood (\ref{eq:ComplexLikelihood}). It is easily seen that this estimator is obtained by solving the equation 
$$\frac{1}{\sqrt{n}} K_n = \psi_1(\sqrt{n} \theta),$$
with $\psi_1(x) = x + \frac{\phi(x)}{\Phi(x)}$, in the case $N_n = n$, and the equation
$$\frac{1}{\sqrt{2n}} K_{2n} = \psi_2(\sqrt{n} \theta),$$
with $\psi_2(x) = x\sqrt{2} + \frac{1}{\sqrt{2}}\frac{\phi(x)}{1 - \Phi(x)}$, in the case $N_n = 2n$. One checks numerically that the map $\psi_1$ strictly increases on $\mathbb{R}$ from $0$ to $\infty$ and that the map $\psi_2$ strictly increases on $\mathbb{R}$ from $-\infty$ to $\infty$. In particular, $\psi_1$ and $\psi_2$ are bijective, from which it follows that $\widehat{\mu}_{c,N_n}$ is uniquely defined by
\begin{equation}
\widehat{\mu}_{c,n} = \frac{1}{\sqrt{n}} \psi_1^{-1}\left(\frac{1}{\sqrt{n}} K_n\right)\label{eq:mucn1}
\end{equation}
if $N_n = n$, and 
\begin{equation}
\widehat{\mu}_{c,2n} =  \frac{1}{\sqrt{n}} \psi_2^{-1}\left(\frac{1}{\sqrt{2n}} K_{2n}\right)\label{eq:mucn2}
\end{equation}
if $N_n = 2n$. Applying the Transformation Theorem, and using (\ref{eq:JointDensity1}) and (\ref{eq:mucn1}), we get, for each $n \in \mathbb{N}_0$ and each $N \in \mathbb{N}_0$,
\begin{equation*}
\mathbb{E}[\left|\widehat{\mu}_{c,N_n}\right|] \geq \frac{1}{\sqrt{n}}\int_{\psi_1(-N)}^{\psi_1(0)} \left|\psi_1^{-1}(u)\right| \phi(u)du \geq - \frac{1}{\sqrt{n}}\phi(\psi_1(0))\int_{\psi_1(-N)}^{\psi_1(0)} \psi_1^{-1}(u)du,\label{eq:MAEmucn}
\end{equation*}
which, by the well known integral equality $\int_a^b f(x) dx + \int_{f(a)}^{f(b)} f^{-1}(x) dx = b f(b) - a f(a),$
$$= \frac{1}{\sqrt{n}}  \phi(\psi_1(0)) \int_{-N}^{0} \psi_1(u) du - N \psi_1(-N),$$
which, plugging in the definition of $\psi_1$ and calculating the integral,
\begin{equation}
= \frac{1}{\sqrt{n}} \phi(\psi_1(0)) \left(\log(1/2) + \frac{N^2}{2} - \log \Phi(-N) - N \frac{\phi(-N)}{\Phi(-N)}\right).\label{eq:BlowUp}
\end{equation}
It can be checked numerically that, for fixed $n$, expression (\ref{eq:BlowUp}) tends to $\infty$ if $N \to \infty$. We conclude that 
\begin{equation}
\forall n: \mathbb{E}[\left|\widehat{\mu}_{c,N_n}\right|] = \infty.\label{eq:mucMAE}
\end{equation}
We infer from (\ref{eq:asymMAEmarginal}) and (\ref{eq:mucMAE}) that conditional bias reduction can cause infinite mean absolute error.


\begin{thebibliography}{99}
\bibitem[BIM18]{BIM18} Berckmoes, B.; Ivanova, A.; Molenberghs, G. (2018) 
{\em On the sample mean after a group sequential trial} 
Computational Statistics \& Data Analysis, 125, 104-118. \url{https://arxiv.org/pdf/1706.01291.pdf}

\bibitem[C89]{C89} 
Chang, M. N. (1989)
{\em Confidence intervals for a normal mean following a group sequential test.}
Biometrics 45, no. 1, 247--254.

\bibitem[EF90]{EF90} 
Emerson, S. S.; Fleming, T. R. (1990)
{\em Parameter estimation following group sequential hypothesis testing.}
Biometrika 77, 875--892.

\bibitem[FDL00]{FDL00}
Fan, X. F.; DeMets, D. L.; Lan, G. (2000)
{\em Bias point of estimation following a group sequential test.}
Technical report \url{https://www.biostat.wisc.edu/sites/default/files/tr_157.pdf}


\bibitem[HP88]{HP88}
Hughes, M.D.; Pocock, S.J. (1988)
{\em Stopping rules and estimation problems in clinical trials.}
Statistics in Medicine 7, 1231--1242.

\bibitem[LH99]{LH99}
Liu, A.; Hall, W. J. (1999)
{\em Unbiased estimation following a group sequential test.}
Biometrika 86, 71--78.

\bibitem[MKA14]{MKA14} Molenberghs, G.; Kenward, M. G.; Aerts, M.; Verbeke, G.; Tsiatis, A. A.; Davidian, M.; Rizopoulos, D.  (2014)
{\em On random sample size, ignorability, ancillarity, completeness, separability, and degeneracy: sequential trials, random sample sizes, and missing data.} 
Stat. Methods Med. Res. 23, no. 1, 11--41.

\bibitem[S78]{S78}
Siegmund, D. (1978) 
{\em Estimation following sequential tests.} 
Biometrika  64, 191--199.

\bibitem[OF79]{OF79}
O'Brien, P. C.; Fleming, T. R. (1979) 
{\em A multiple testing procedure for clinical trials.} 
Biometrics 35, 549-556.

\bibitem[P77]{P77}
Pocock, S. J. (1977) 
{\em Group sequential methods in the design and analysis of clinical trials.} 
Biometrika. 64 (2): 191-9.

\bibitem[W92]{W92}
Woodroofe, M. (1992)
{\em Estimation after sequential testing: a simple approach for a truncated sequential probability ratio test.} 
Biometrika 79, no. 2, 347--353.
\end{thebibliography}
\end{document}